\documentclass[11pt]{paper}
\usepackage{amssymb,euscript,graphicx,pifont,cite,amsmath,epsfig}

\usepackage{cite}
\usepackage{pifont}
\usepackage{times}
\usepackage{amssymb,epsfig,floatflt}
\usepackage{subfigure}[1995/03/06]
\usepackage[section]{placeins}
\usepackage{epsfig,psfrag}
\usepackage{graphics}
\usepackage{color}

\input amssymb.sty
\newcommand{\II}{\mathbb{I}}

\def\dspace{\multiply\normalbaselineskip 150

          \divide\normalbaselineskip 100 \normalbaselines

\csname @@normalbaselineskip\endcsname\normalbaselineskip}

\def\rr{\rm I\!R}

\newtheorem{thm}{Theorem}

\newtheorem{lem}{Lemma}

\newtheorem{defn}{Definition}

\newtheorem{cor}{Corollary}

\usepackage{cite}

\begin{document}
\thispagestyle{empty}

\title{\LARGE \bf Stability Margins of $\mathcal{L}_1$ Adaptive Controller: Part II\thanks{Research is supported by AFOSR under Contract
No. FA9550-05-1-0157.}}
\author{Chengyu Cao and Naira Hovakimyan\thanks{The authors are with
  Aerospace \& Ocean Engineering, Virginia Polytechnic
 Institute \& State University, Blacksburg, VA 24061-0203, e-mail:
  chengyu, nhovakim@vt.edu}}
\date{}
\maketitle
\begin{abstract}
In Part I of this paper, \cite{CDC06_chengyu_L1}, we have
developed a novel $\mathcal{L}_1$ adaptive control architecture
that enables fast adaptation and leads to uniformly bounded
transient and asymptotic tracking for system's both signals, input
and output, simultaneously. In this paper, we derive the stability
margins of $\mathcal{L}_1$ adaptive control architecture,
including time-delay and gain margins in the presence of
time-varying bounded disturbance.
 Simulations
 verify the theoretical findings.
\end{abstract}

\section{Introduction}

Adaptive control  schemes have proven to be extremely useful in a
number of flight tests for recovering the nominal performance in
the presence of  modeling and environmental uncertainties (see
\cite{wise} and references therein). A major challenge in analysis
of these systems is
 determining its stability margins dependent upon the adaptation
gain. Today it largely relies on the numerical evidence provided
by Monte-Carlo schemes. It has been observed that  increasing the
adaptation gain leads to improved tracking performance, but
results in high-frequency oscillations in the control signal and
reduces the system's tolerance  to the time-delay in the control
and the sensor channels.
%
%
%

In the linear time invariant (LTI) systems theory, stability
margins are defined by the gain and the  phase  margins. Phase
margin characterizes the amount of additional phase lag at the
gain-crossover frequency required to bring the system to the verge
of instability. Phase margin is significant in predicting how much
time-delay the system can endure in its input/output channels
before it loses its stability. While the gain margin can be
generalized for nonlinear systems, the notion of the phase margin
cannot be extended to nonlinear systems in straightforward manner.
Instead it is common to use sector and disk margins for nonlinear
systems \cite{sepulchre97}. However, from practical control design
perspective these notions are not as useful as the phase margin in
the linear systems theory. In this paper, instead of the phase
margin we introduce the notion of the time-delay  margin directly
for the closed-loop nonlinear adaptive system. Time-delay margin
characterizes  the maximum time-delay in the (sensor) channel that
the closed-loop system can tolerate before it loses its stability.
In linear systems theory this corresponds to the ratio of the
phase margin to the cross-over frequency of its  Bode plot.
Similarly, the gain margin is the maximum open loop gain that the
closed-loop system can tolerate before it loses its stability.

In \cite{CaoHovakimyan, CaoHovakimyan1}, we have introduced novel
$\mathcal L_1$ adaptive control architecture that has guaranteed
transient performance. In \cite{CDC06_chengyu_L1}, we have
extended the approach to systems with unknown time-varying
parameters and bounded disturbances. In this paper, we derive the
stability margins for the $\mathcal L_1$ adaptive control
architecture from \cite{CDC06_chengyu_L1}, which we specialize for
unknown  constant  parameters and bounded time-varying
disturbances. While the analysis of the gain-margin is relatively
straightforward, the analysis of its time-delay margin takes
several steps. At first we introduce an equivalent linear-time
invariant (LTI) system, subject to an exogenous input dependent
upon the parameters and time trajectories of certain signals of
the closed-loop adaptive system. We prove that with the same
initial conditions in the presence of the same time-delay in the
output channels of these two systems there exists at least one
exogenous input  such that their resulting trajectories are the
same over the entire time-horizon. Next, we prove that for every
value of the time-delay within the time-delay margin of this LTI
system there exists a lower bound for the adaptive gain that
renders this exogenous input bounded.

 We notice that  characterization
of the time-delay margin is extremely difficult as compared to the
gain-margin analysis for nonlinear closed-loop systems. To the
best of our knowledge there are no such results in adaptive
control theory, despite the fact that there is a large body of
well-established literature on adaptive control of time-delay
systems. Control of time-delay systems and determining the
time-delay margin of a closed-loop system are principally
different problems, and one cannot be used to provide a solution
for the other. On the other hand, this is not surprising since the
time-delay margin cannot be characterized if the transient is not
guaranteed.

The paper is organized as follows. Section \ref{sec:preliminary}
states some preliminary definitions, and Section \ref{sec:PF}
gives the problem formulation. In Section \ref{sec:Hinf}, the
$\mathcal{L}_1$ adaptive controller  is presented. Stability
margins, including time-delay and gain margins, are derived in
Section \ref{sec:convergence}. Results of \cite{CDC06_chengyu_L1}
and of this paper are generalized in Section
\ref{sec:generalization}. In section \ref{sec:simu}, simulation
results are presented, while Section \ref{sec:conclusion}
concludes the paper. The proof of the main theorem is in Appendix.

\section{Preliminaries}\label{sec:preliminary}

In this Section, we recall some basic definitions and facts from
linear systems theory, \cite{IoaBook03,KhaBook02,ZhoBook98}.

%
%
%
%
%
%
%
%
%
%

\begin{defn}\label{defn1}
For a signal $\xi(t)=[\xi_1(t),\cdots,\xi_n(t)]^\top,~ t\geq 0$,
its truncated ${\mathcal L}_{\infty}$ norm and ${\mathcal
L}_{\infty}$ norm are defined as
$
 \Vert \xi_t \Vert_{{\mathcal L}_\infty}  =
\max_{i=1,..,n} (\sup_{0\leq \tau \leq t} |\xi_i(\tau)|)$, $ \Vert
\xi \Vert_{{\mathcal L}_\infty}  =  \max_{i=1,..,n}
(\sup_{\tau\geq 0} |\xi_i(\tau)|)$.
\end{defn}
\begin{defn}
 The $\mathcal{L}_1$ gain of a stable proper single--input single--output system $H(s)$
 is defined as
$|| H(s)||_{\mathcal{L}_1} = \int_{0}^{\infty} |h(t)| d t, $ where
$h(t)$ is the impulse response of $H(s)$, computed via the inverse
Laplace transform $ h(t)=\frac{1}{2 \pi i}\int_{\alpha-i
\infty}^{\alpha+i \infty} H(s) e^{st} ds, t\geq 0, $ in which
integration is done along the vertical line $x=\alpha>0$ in
complex plane.
\end{defn}

{\em Proposition:}
 A continuous time LTI
system (proper) with impulse response $h(t)$ is stable if and only
if
$
\int_{0}^{\infty} |h(\tau)| d\tau < \infty.
$
A proof can be found in \cite{IoaBook03} (page 81, Theorem 3.3.2).

\begin{defn}
For a stable proper $m$ input $n$ output system $H(s)$ its
$\mathcal{L}_1$ gain is defined as
$
\Vert H(s) \Vert_{\mathcal{L}_1} = \max_{i=1,..,n} \left(
\sum_{j=1}^{m}\Vert H_{ij}(s)\Vert_{\mathcal{L}_1} \right)\,,
$
where $H_{ij}(s)$ is the $i^{th}$ row $j^{th}$ column element of
$H(s)$.
\end{defn}


\begin{lem}\label{lem:L1} For a stable proper multi-input multi-output (MIMO)  system $H(s)$ with input $r(t) \in {\rr}^m$ and output
 $x(t)\in {\rr}^n$, we have
$ \Vert x_t\Vert_{{\mathcal L}_\infty} \leq \Vert
H\Vert_{{\mathcal L}_1} \Vert r_t\Vert_{{\mathcal L}_\infty},
\forall t\geq 0$.
\end{lem}
%
%

\begin{cor}\label{lem:L1_ext} For a stable proper MIMO  system  $H(s)$, if the input   $r(t) \in {\rr}^m$ is bounded, then the  output
 $x(t)\in {\rr}^n$ is also bounded as
$ 
\Vert x\Vert_{{\mathcal L}_\infty } \leq \Vert
H(s)\Vert_{\mathcal{L}_1} \Vert r\Vert_{{\mathcal L}_\infty }.
$ 
\end{cor}

Consider a linear time invariant system: $ \dot{x}(t) = A x(t) +b
u(t) $, where $x\in {\rr}^n$, $u\in {\rr}$, $b\in {\rr}^{n}$,
$A\in {\rr}^{n \times n}$ is Hurwitz, and assume  $(s I -A)^{-1} b
$ is strictly proper and stable. Notice that it can be expressed
as: $ (s I -A)^{-1} b = \frac{n(s)}{d(s)} $, where
$d(s)={\rm{det}} (s I- A)$ is a $n^{th}$ order stable polynomial,
and $n(s)$ is a $n\times 1$ vector with its $i^{th}$ element being
a polynomial function: $ n_i(s) = \sum_{j=1}^{n} n_{ij} s^{j-1} $.

\begin{lem}\label{lem:pre1}
If $(A\in {\rr}^{n\times n}, b\in {\rr}^{n})$ is controllable, the
matrix $N$ with its $i^{th}$ row $j^{th}$ column entry $n_{ij}$ is
full rank.
\end{lem}

\begin{lem}\label{lem:pre2}
If $(A, b)$ is controllable and $(s I -A)^{-1} b $ is strictly
proper and stable, there exists $c\in {\rr}^n$ such that
$c^{\top}(s I -A)^{-1} b$ is minimum phase
 with relative degree one, i.e. all its zeros are located in
the left half plane, and its denominator is one order larger than
its numerator.
\end{lem}

Also, we introduce the following notations that will be used
throughout the paper. Let $x_h(t)$ be the state variable of the
LTI system $H_x(s)$, while $x_i(t)$ and $x_s(t)$ be the input and
the output signals of it. We note that for any time instant $t_1$
and any fixed time-interval $[t_1,\; t_2]$, where $t_2>t_1$, given
$x_h(t_1)$ and an impulse-free input signal $x_i(t)$ over
$[t_1,\,t_2)$, $x_s(t)$ is uniquely defined for  $t\in [t_1,
t_2]$. Let $ \mathcal{S}$ be the map $ x_s(t)|_{t\in[t_1,\;t_2]} =
\mathcal{S} ( H_x(s), x_h(t_1), x_i(t)|_{t\in[t_1,\;t_2)} )$. We
note that $\mathcal{S}$ is continuous, if $x_i(t)$ is impulse
free. Also, $x_s(t)$ is defined over a closed interval
$[t_1,\;t_2]$, although $x_i(t)$ is defined over the corresponding
open set $[t_1,\;t_2)$. The next lemma follows from the definition
of $\mathcal S$.
\begin{lem}\label{Smaplem}
Let $ x_{o_1}|_{t\in[t_1,\;t_2]}  =  \mathcal{S} ( H_x(s),
x_{h_1}, x_{i_1}(t)|_{t\in[t_1,\;t_2)} )$, $
x_{o_2}|_{t\in[t_1,\;t_2]}  =  \mathcal{S} ( H_x(s), x_{h_2},
x_{i_2}(t)|_{t\in[t_1,\;t_2)} )$. If $ x_{h_1} = x_{h_2}$ and
$x_{i_1}(t)=x_{i_2}(t) $ over $[t_1\,,t_2)$, then $
x_{o_1}(t)=x_{o_2}(t)$ for any $ t\in[t_1\,,t_2]$.
\end{lem}

\section{Problem Formulation}\label{sec:PF}
Consider the following single-input single-output system dynamics:
\begin{eqnarray}
\dot{x}(t)  & = &  A_m x(t)+b \left( \omega u(t)+\theta^{\top}
x(t)+ \sigma(t) \right),
x(0)=x_0\,\nonumber\\
y(t) & = & c^{\top} x(t)\,,\label{problemnew}
\end{eqnarray}
where $x\in {\rr}^n$ is the system state vector (measurable),
$u\in {\rr}$ is control signal, $y\in {\rr}$ is the regulated
output, $b,c\in {\rr}^n$ are known constant vectors, $A_m\in
\rr^{n \times n}$ is given Hurwitz matrix,  $\omega\in \rr$ is
unknown constant with given sign, $\theta \in \rr^n$ is unknown
constant vector, and $\sigma(t)\in \rr$ is a uniformly bounded
time-varying disturbance with a uniformly bounded derivative.
Without loss of generality, we assume
\begin{equation}\label{Thetadef}
\omega\in \Omega_0=[\omega_{l_0},\;\omega_{u_0}]\,,~\theta \in
\Theta\,,~|\sigma(t)| \leq \Delta_0\,, ~ \forall ~t\geq 0\,,
\end{equation}
where $\omega_{u_0}>\omega_{l_0}>0$ are known (conservative) upper
and lower bounds, $\Theta$ is a known (conservative) compact set
and $\Delta_0\in \rr^{+}$ is a known (conservative)
$\mathcal{L}_{\infty}$ bound of $\sigma(t)$. We further assume
that $\sigma(t)$ is continuously differentiable and its derivative
is uniformly bounded, i.e. $ | \dot{\sigma}(t) |  \leq
d_{\sigma}<\infty$ for any $t\geq 0$, where $d_{\sigma}$ can be
arbitrarily large as long as it is finite.

In \cite{CDC06_chengyu_L1}, we have considered the system in
(\ref{problemnew}) in the presence of time-varying  $\theta(t)$
and have designed an adaptive controller to ensure that $y(t)$
tracks a given bounded continuous reference signal $r(t)$ {\em
both in transient and steady state}, while all other error signals
remain bounded. The main result of \cite{CDC06_chengyu_L1} implies
that by increasing the adaptation gain one can get arbitrarily
close transient and asymptotic tracking of a desired reference
system. In \cite{CDC06_chengyu_L1}, we have also considered the
particular case of constant $\theta$. Here we investigate the
stability margins for this latter case.

\section{ $\mathcal{L}_1$ Adaptive Controller}\label{sec:Hinf}

In this section, we repeat  the $\mathcal{L}_1$ adaptive control
architecture for the system in (\ref{problemnew}) that permits
complete transient characterization for both $u(t)$ and $x(t)$,
 \cite{CDC06_chengyu_L1}.
The elements of $\mathcal{L}_1$ adaptive controller are introduced
next without repeating the  proofs from \cite{CDC06_chengyu_L1}.

{\bf Companion Model:} The companion model is:
\begin{eqnarray}
\dot{\hat{x}}(t) & = & A_m \hat{x}(t)+b ( \hat{\omega}(t) u(t) +
\hat{\theta}^{\top}(t) x(t)+\hat{\sigma}(t) )\,,
\nonumber\\
\hat{y}(t) & = & c^{\top} \hat{x}(t)\,,\quad \hat
x(0)=x_0\,,\label{L1_companionmodal}
\end{eqnarray}
which has the same dynamic structure as the system in
(\ref{problemnew}). Only the unknown parameters and the
disturbance $\omega, \theta, \sigma(t)$ are replaced by their
adaptive estimates $\hat{\omega}(t), \hat{\theta}(t),
\hat{\sigma}(t)$.

{\bf Adaptive Laws:} Adaptive estimates are governed by the
following  laws:
\begin{eqnarray}
\dot{\hat{\theta}}(t) &  = & -\Gamma_{\theta}{\rm Proj}(x(t)
\tilde{x}^{\top}(t) P b , \hat{\theta}(t)), \label{adaptivelaw_L1}\\
\dot{\hat{\sigma}}(t) &  = & -\Gamma_{\sigma}{\rm Proj}(
\tilde{x}^{\top}(t) P b , \hat{\sigma}(t)),  \label{adaptivelaw_L2}\\
\dot{\hat{\omega}}(t) &  = & -\Gamma_{\omega}{\rm Proj}( u(t)
\tilde{x}^{\top}(t) P b, \hat{\omega}(t)), \label{adaptivelaw_L3}
\end{eqnarray}
where $\tilde{x}(t)=\hat x(t)-x(t)$ is the error signal between
the state of the system and the companion model, $P$ is the
solution of the algebraic equation $ A_m^{\top} P+P A_m =- Q$,
$Q>0$,  $\Gamma_{\theta}=\Gamma_c \II_{n\times n} \in
{\rr}^{n\times n}$, $\Gamma_{\sigma}=\Gamma_{\omega}=\Gamma_c$ are
adaptation gains with $\Gamma_c\in \rr^{+}$. In the implementation
of the projection operator  we use the compact sets $\Theta$  as
given in (\ref{Thetadef}), while we replace  $\Delta_0$,
$\Omega_0$ by larger sets $\Delta$ and
$\Omega=[\omega_{l},\;\omega_{u}]$ such that
\begin{equation}\label{eqn:Delta}
\Delta_0< \Delta, ~~\, 0<\omega_l<\omega_{l_0}<\omega_{u_0}<
\omega_u \,.
\end{equation}

 The purpose of this will be shortly clarified in the analysis of the stability margins.

{\bf Control Law:} The control signal is generated through gain
feedback of the following  system:
\begin{eqnarray}
\chi(s)  =  D(s) r_u(s) \,, \quad u(s)  =  - k
\chi(s)\,,\label{controllaw}
\end{eqnarray}
where $r_u(s)$ is the Laplace transformation of
$r_u(t)=\hat{\omega}(t) u(t)+\bar{r}(t)$, $k\in \rr^{+}$ is a
feedback gain, $ \bar{r}(t)=\hat{\theta}^{\top}(t)
x(t)+\hat{\sigma}(t)-k_g r(t)$, $
 k_g =-\frac{1}{ c^{\top} A_m^{-1} b}$,
and $D(s)$ is a  LTI system that needs to be chosen  to ensure
\begin{equation}\label{Csdef}
 C(s) = \frac{\omega k D(s)}{1+\omega k D(s)}
\end{equation}
is stable and strictly proper with  $C(0)=1$. One  choice is $
D(s) =\frac{1}{s}$, that leads to $ C(s) = \frac{\omega k
}{s+\omega k }$. Let $ L = \max_{\theta\in \Theta} \sum_{i=1}^{n}
|\theta_i|$. We now give the $\mathcal{L}_1$ performance
requirement that ensures
 desired transient performance,  \cite{CDC06_chengyu_L1}.

{\bf $\mathcal{L}_1$-gain stability requirement:} Design $D(s)$ to
ensure that $C(s)$ in (\ref{Csdef}) satisfies
\begin{equation}\label{condition3}
\Vert G(s) \Vert_{\mathcal{L}_1} L <1\,,
\end{equation}
where $G(s)=H(s)(1-C(s))$, and
$
H(s) = (s I -A_m)^{-1} b\,.
$

The complete $\mathcal{L}_1$ adaptive controller consists of
(\ref{L1_companionmodal}),
(\ref{adaptivelaw_L1})-(\ref{adaptivelaw_L3}), (\ref{controllaw})
subject to (\ref{condition3}).
We notice that the $\mathcal{L}_1$-gain stability requirement
depends only upon the choice of $\Theta$ and is independent of the
choice of  $\Delta_0$, $\Omega_0$ or $\Delta$, $\Omega$.

\section{Analysis of $\mathcal{L}_1$ Adaptive
Controller}\label{sec:convergence}

 Next, consider the following
closed-loop reference system with the control signal $u_{ref}(t)$
and the system response $x_{ref}(t)$, the stability of which,
subject to (\ref{condition3}), can be proven using the small-gain
theorem, \cite{CDC06_chengyu_L1}:
\begin{eqnarray}
\dot{x}_{ref}(t) & =&  A_m x_{ref}(t)+b ( \omega u_{ref}(t)+\theta^{\top} x_{ref}(t)+\sigma(t))\nonumber\\
 u_{ref}(s)  &=&  C(s) \frac{\bar{r}_{ref}(s)}{\omega}\,, \quad
y_{ref}(t)=c^\top x_{ref}(t)\,, \quad\,\,  \label{refy}
\end{eqnarray}
with $x_{ref}(0)=x_0$, where  $\bar{r}_{ref}(s)$ is the Laplace
transformation of the signal $ \bar{r}_{ref} = -\theta^{\top}
x_{ref}(t)-\sigma(t)+ k_g r(t)$.

\begin{lem}\label{lem:1}\cite{CDC06_chengyu_L1}
For the system in (\ref{problemnew}) and the $\mathcal{L}_1$
adaptive controller in (\ref{L1_companionmodal}),
(\ref{adaptivelaw_L1})-(\ref{adaptivelaw_L3})  and
(\ref{controllaw}), we have $ \Vert \tilde{x}
\Vert_{\mathcal{L}_{\infty}}  \leq
\sqrt{\frac{\theta_{m}}{\lambda_{\min}(P) \Gamma_c}}$ where $
\theta_{m} \triangleq \max_{\theta\in \Theta} \sum_{i=1}^{n} 4
\theta_i^2+4\Delta^2+4 \left(\omega_{u}-\omega_{l} \right)^2+ 2
\frac{\lambda_{\max}(P)}{\lambda_{\min}(Q)}  b_{\sigma} \Delta $ .
\end{lem}

Lemma \ref{lem:pre2} ensures existence of $c_{o}\in {\rr}^n$ such
that $ c_{o}^{\top} H(s) = \frac{N_n(s)}{N_d(s)}$, where the order
of $N_d(s)$ is one more than the order of $N_n(s)$, and both
$N_n(s)$ and $N_d(s)$ are stable polynomials.
\begin{thm}\label{thm:5}\cite{CDC06_chengyu_L1}
Given the system in (\ref{problemnew}) and the $\mathcal{L}_1$
adaptive controller defined via (\ref{L1_companionmodal}),
(\ref{adaptivelaw_L1})-(\ref{adaptivelaw_L3}) and
(\ref{controllaw}) subject to (\ref{condition3}), we have: $ \Vert
x-x_{ref} \Vert_{{\mathcal L}_{\infty}}   \leq  \gamma_1, \quad
\Vert u - u_{ref} \Vert_{{\mathcal L}_{\infty}}  \leq \gamma_2$,
where $ \gamma_1  =  \frac{\Vert C(s)
\Vert_{\mathcal{L}_1}}{1-\Vert G(s)\Vert_{\mathcal{L}_1} L}
\sqrt{\frac{\theta_{m}}{\lambda_{\max}(P) \Gamma_c}}$, and $
\gamma_2  =  \Vert \frac{C(s)}{\omega} \theta^{\top}
\Vert_{\mathcal{L}_1} \gamma_1 + \Big\| \frac{C(s)}{\omega}
\frac{1}{c_{o}^{\top} H(s)} c_{o}^{\top} \Big\|_{\mathcal{L}_1}
\sqrt{\frac{\theta_{m}}{\lambda_{\max}(P) \Gamma_c}}$.
\end{thm}

\section{Time-delay Margin Analysis}

\subsection{$\mathcal{L}_1$ adaptive controller in the presence of
 time-delay}

To analyze the time-delay margin of  the closed-loop adaptive
system in the next section we consider a linear time-invariant
(LTI) system subject to an external exogenous input. We develop
sufficient conditions under which that LTI system can be used to
evaluate the time-delay margin of the closed-loop adaptive system.
Before then, we need to introduce the following three systems.

{\bf System 1.} We rewrite the open-loop system in
(\ref{problemnew}) as
\begin{equation}\label{L1_td1}
x(s) = \bar{H}(s) ( \omega u(s)+ \sigma(s))\,,
\end{equation}
where $ \bar{H}(s)  =  (s I- A_m- b \theta^{\top})^{-1} b$.
Without loss of generality, we set:
\begin{equation}
x(0)=0.\label{zeroinit0}
\end{equation}
Let $x_d(t)$ be the  delayed signal of the open-loop state $x(t)$
by a constant time interval $\tau$, i.e
\begin{equation}
x_d(t) =\left\{
\begin{array}{ll}
x(t-\tau) & t\geq \tau\,, \\
0 & t<\tau\,.
\end{array}
\right. \label{zeroinit1}
\end{equation}
We close the loop of (\ref{L1_td1}) with
 $\mathcal L_1$ adaptive controller
(\ref{L1_companionmodal}),
(\ref{adaptivelaw_L1})-(\ref{adaptivelaw_L3}), (\ref{controllaw}),
using $x_d(t)$ from (\ref{zeroinit1}) instead of $x(t)$ everywhere
in the definition of (\ref{L1_companionmodal}),
(\ref{adaptivelaw_L1})-(\ref{adaptivelaw_L3}), (\ref{controllaw}).
We denote the resulting control and state trajectories of this
closed-loop system by $u(t)$ and $x_{d}(t)$. We further notice
that this closed-loop adaptive system has a unique solution. It is
the stability of this closed-loop system that we are investigating
in this paper, dependent upon $\tau$. It is important to point out
that while applying the $\mathcal L_1$ adaptive controller
(\ref{L1_companionmodal}),
(\ref{adaptivelaw_L1})-(\ref{adaptivelaw_L3}), (\ref{controllaw})
to the system in (\ref{L1_td1}) using $x_d(t)$ from
(\ref{zeroinit1}), one cannot derive the dynamics of the error
signal between the system state and the companion model, the
boundedness of which is stated in Lemma \ref{lem:1}. Neither
Theorem \ref{thm:5} is valid.

{\bf System 2.} Next, we consider the following closed-loop
system with the same zero initial conditions:
\begin{equation}
\dot{x}_q(t)  =  A_m x_q(t) + b \left(\omega u_q(t)+\theta^{\top}
x_q(t) + \sigma(t)+\eta(t)  \right)\,,\label{L1_tv}
\end{equation}
where $x_q(0)=x(0)$, $\theta$ and $\sigma(t)$  have been
introduced in (\ref{problemnew}), $u_q(t)$ is defined via
(\ref{L1_companionmodal}),
(\ref{adaptivelaw_L1})-(\ref{adaptivelaw_L3}) and
(\ref{controllaw}), while $\eta(t)$ is a continuously
differentiable bounded signal with uniformly bounded derivative.
As compared to (\ref{problemnew}) or (\ref{L1_td1}), the system in
(\ref{L1_tv}) has one more additional disturbance signal
$\eta(t)$. If
\begin{eqnarray}\label{disinset}
|\sigma(t)+\eta(t)|  \leq  \Delta\,,
\end{eqnarray}
where $\Delta$ has been defined in (\ref{eqn:Delta}), then
application of $\mathcal L_1$ adaptive controller to the system in
(\ref{L1_tv}) is well defined, and hence the results of Theorem
\ref{thm:5} are valid for the system in (\ref{L1_tv}) as well. We
denote by $u_q(t)$ the time trajectory of the $\mathcal L_1$
adaptive controller, resulting from its application to
(\ref{L1_tv}).

{\bf System 3.} Finally, we consider the open-loop system in
(\ref{L1_td1})-(\ref{zeroinit1}) and apply $u_q(t)$  to it and
look at its delayed output $x_{o}(t)$, where the subindex $o$ is
added to indicate the open-loop nature of this signal. It is
important to notice that at this point we view $u_q(t)$ as a
time-varying input signal for (\ref{L1_td1}), and not as a
feedback signal, so that (\ref{L1_td1}) remains an open-loop
system in this context.

\begin{figure}[!h]
\begin{center}
\includegraphics[width=3.0in,height=2.0in]{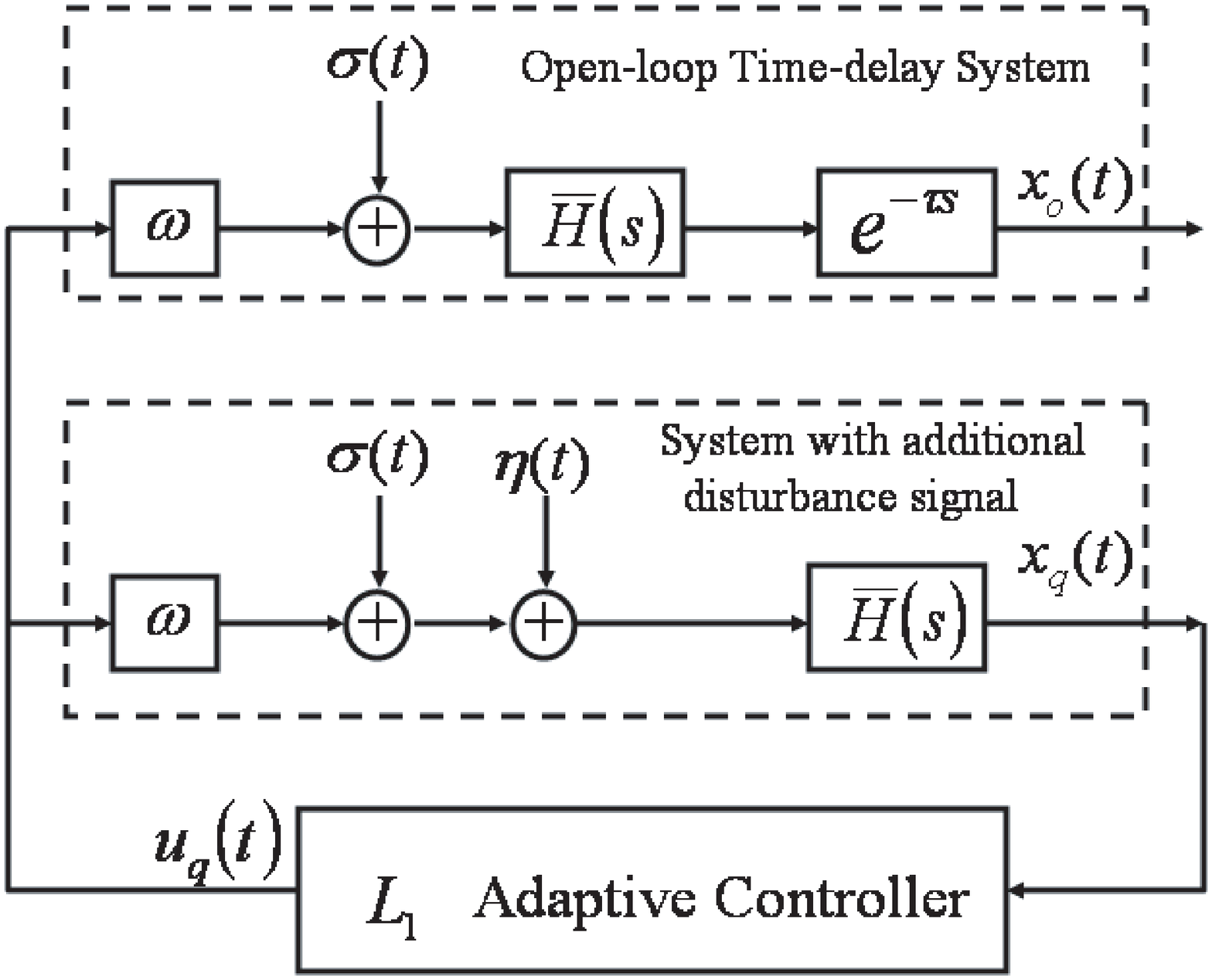}
\caption{Systems 2 and 3}\label{fig:delayeq}
\end{center}
\end{figure}

Illustration of these last two systems is given
 in Fig. \ref{fig:delayeq}.

\begin{lem}\label{lem:equivalence}
 If the time-delayed output of the open-loop System 3 has the same
 time history as the closed-loop output of System 2, i.e.
\begin{equation}\label{output_eq}
x_{o}(t)=x_q(t), \;\forall ~t\geq 0\,,
\end{equation}
then $ u(t)=u_q(t),\;x_{d}(t)=x_q(t)$, ~ $\forall ~t\geq 0$, where
$u(t)$ and $x_{d}(t)$ denote the control and state trajectories of
the closed-loop System 1 in (\ref{L1_td1})-(\ref{zeroinit1}) with
$\mathcal L_1$ adaptive controller.
\end{lem}
{\bf Proof.} It follows from (\ref{output_eq}) that the open-loop
time-delayed System  3 in (\ref{L1_td1})-(\ref{zeroinit1})
generates $x_q(t)$ in response to the input $u_q(t)$. When applied
to (\ref{L1_tv}), $u_q(t)$ leads to $x_q(t)$. Hence,
 $u_q(t)$ and $x_q(t)$
 are also  solutions of the closed-loop adaptive System 1 in
(\ref{L1_td1})-(\ref{zeroinit1}) with (\ref{L1_companionmodal}),
(\ref{adaptivelaw_L1})-(\ref{adaptivelaw_L3}), (\ref{controllaw}).
$\hfill{\square}$

This Lemma consequently implies that to ensure stability of the
System 1 in the presence of a given  time-delay $\tau$, it is
sufficient to prove existence of $\eta(t)$ in System 2, satisfying
(\ref{disinset}) and verifying (\ref{output_eq}). We notice,
however, that the closed-loop System 2 is a nonlinear system due
to the nonlinear adaptive laws, so that the proof on existence of
such $\eta(t)$ for this system and explicit construction of the
set $\Delta$ is not straightforward. Moreover, we note that the
condition in (\ref{output_eq}) relates the time-delay $\tau$ of
System 1 (or System 3) to the  signal $\eta(t)$ implicitly. In the
next section of this paper we introduce an equivalent  LTI system
that helps to prove existence of such $\eta(t)$ and leads to
explicit construction of $\Delta$. Definition of this LTI system
is the key step in the overall analysis. It has an exogenous input
that lumps the time trajectories of  the nonlinear elements of the
closed-loop System 2.
 For this LTI system, the time delay margin can be
computed via its open-loop transfer function, which consequently
defines a conservative lower bound for the time-delay margin of
the adaptive system.

\subsection{LTI System in the Presence of Time-delay in its Output}
Consider the following
closed-loop LTI system:
\begin{eqnarray}
& & x_l(s)  =  \bar{H}(s)
\zeta_{l}(s),\;\;\epsilon_l(s)  =  (C(s)/\omega) \tilde{r}_l(s)\nonumber\\
& & u_l(s)  = (1/\omega) C(s) ( k_g r(s)-  \theta^{\top}
x_l(s)-\sigma(s)-\eta_l(s) )-\epsilon_l(s) \nonumber
\end{eqnarray}
where $\zeta_l(s)  =  \omega u_l(s) +\sigma(s)$, $ \eta_l(s)  =
\zeta_{l}(s)-\omega u_l(s) -\sigma(s)$,
$r(s)$ and $\sigma(s)$ are the Laplace transformations of the
bounded signals $r(t)$ and $\sigma(t)$, respectively,  $x_l(t)$,
$u_l(t)$ and $\epsilon_l(t)$ are selected states, $\zeta_l(t)$ is
its output signal, and $\tilde{r}_l(s)$ is the Laplace
transformation of an exogenous
 signal $\tilde{r}_l(t)$.  We note that the
system trajectories are uniquely defined once $\tilde r_l(t)$ is
given.

We notice that  the transfer functions from $\sigma(t)$ and $r(t)$
to $x_l(t)$ are the same as in the reference system. Since
$x_l(s)= \bar{H}(s) \zeta_l(s)$, we have
\begin{eqnarray}
x_l(s)/r(s) =   ( \bar{H}(s) C(s))/(1+C(s)\theta^{\top} \bar{H}(s))\,,&& \label{LTI_r}\\
x_l(s)/\sigma(s)  =  (\bar{H}(s)(1-C(s)))/(1+C(s)\theta^{\top}
\bar{H}(s))\,. &&\label{LTI_sigma}
\end{eqnarray}
One can verify  that for the reference system in (\ref{refy}), we
have $x_{ref}(s)/r(s)$ and $x_{ref}(s)/\sigma(s)$ equivalent to
(\ref{LTI_r}) and (\ref{LTI_sigma}). We also notice that the LTI
system without time-delay ensures stable transfer functions from
inputs $r(t)$, $\sigma(t)$ and $\tilde{r}_l(t)$ to output
$\zeta_l(t)$.

\begin{figure}[!h]
\begin{center}
\includegraphics[width=3.0in,height=2.0in]{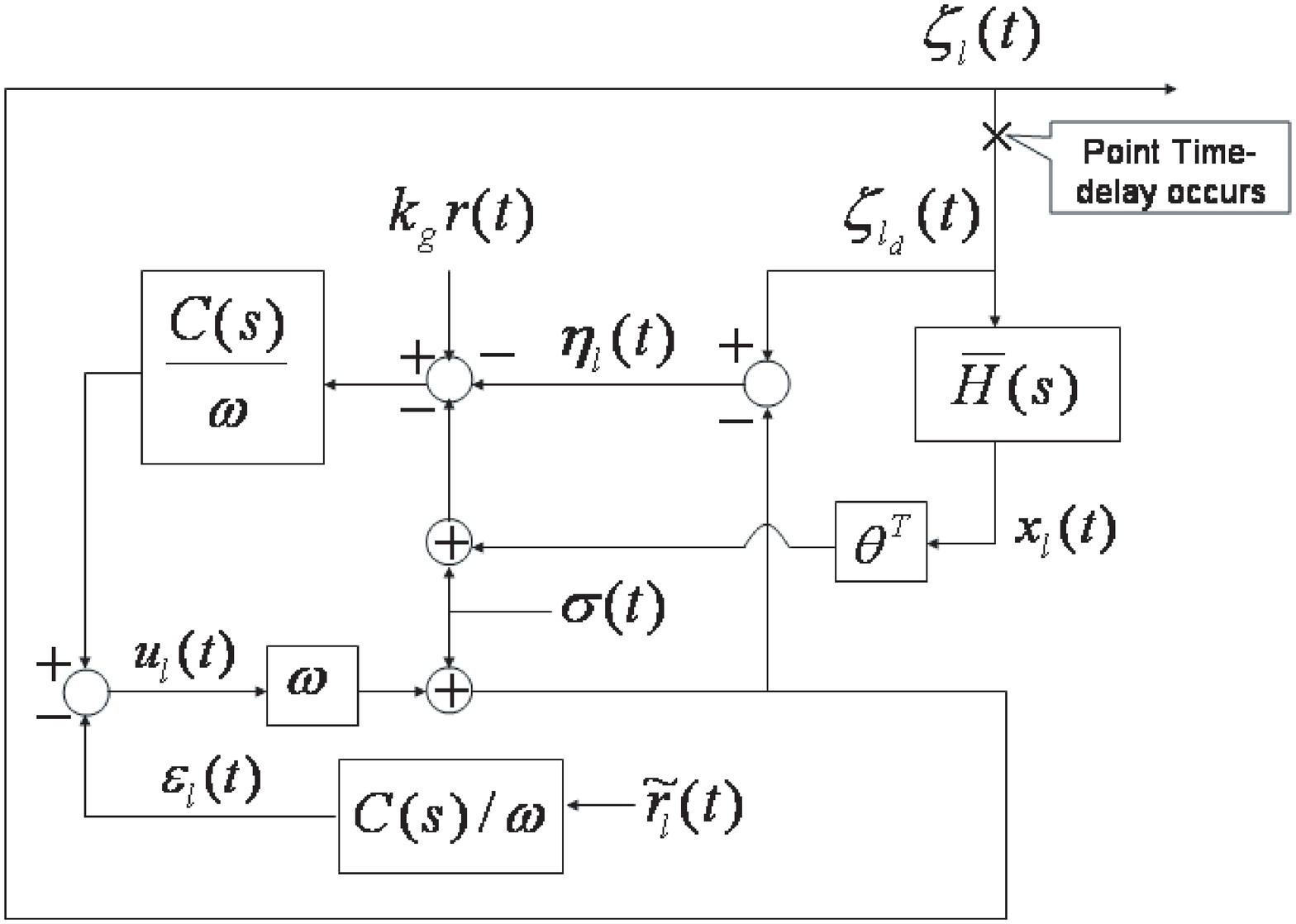}
\caption{LTI system}\label{fig:LTImargin}
\end{center}
\end{figure}

Assume the system output $\zeta_l(t)$ experiences time-delay
$\tau$, so that in the presence of the time-delay we have:
\begin{eqnarray}
&& x_l(s)  =  \bar{H}(s)
\zeta_{l_d}(s)\label{LTImarginori3}\\
&& u_l(s)  =  \scriptstyle{(C(s)/\omega) \left( k_g r(s)-
\theta^{\top} x_l(s)-\sigma(s)-\eta_l(s) \right)-\epsilon_l(s)}
\label{LTImarginori2}\\
&&\epsilon_l(s)  =  (C(s)/\omega) \tilde{r}_l(s) \label{LTImarginori5}\\
 &&\zeta_l(s)  =  \omega u_l(s) +\sigma(s)\,,\label{LTImarginori4}
\end{eqnarray}
where $\zeta_{l_d}(t)$ is the time-delayed
 signal of $\zeta_l(t)$, i.e
\begin{equation}
\zeta_{l_d}(t) =\left\{
\begin{array}{ll}
0 & t<\tau\,, \\
\zeta_l(t-\tau) & t\geq \tau\,,
\end{array}
\right.\label{zeroinit4}
\end{equation}
consequently leading to redefined $\eta_l(s)$:
\begin{equation}\label{LTImarginori1}
 \eta_l(s)  =  \zeta_{l_d}(s)-\omega u_l(s)
-\sigma(s).
\end{equation}
 Let
\begin{equation}
x_l(0)=0,\,\quad  u_l(0)=0\,, \quad
\epsilon_l(0)=0\,.\label{zeroinit3new}
\end{equation}
We notice that the system in
(\ref{LTImarginori3})-(\ref{LTImarginori4}) is highly coupled. Its
diagram is plotted in Figure \ref{fig:LTImargin}.

\subsection{Time-Delay Margin of the LTI System}\label{sec:timedelaymargin}

\begin{figure}[!h]
\begin{center}
\includegraphics[width=2.5in,height=1.4in]{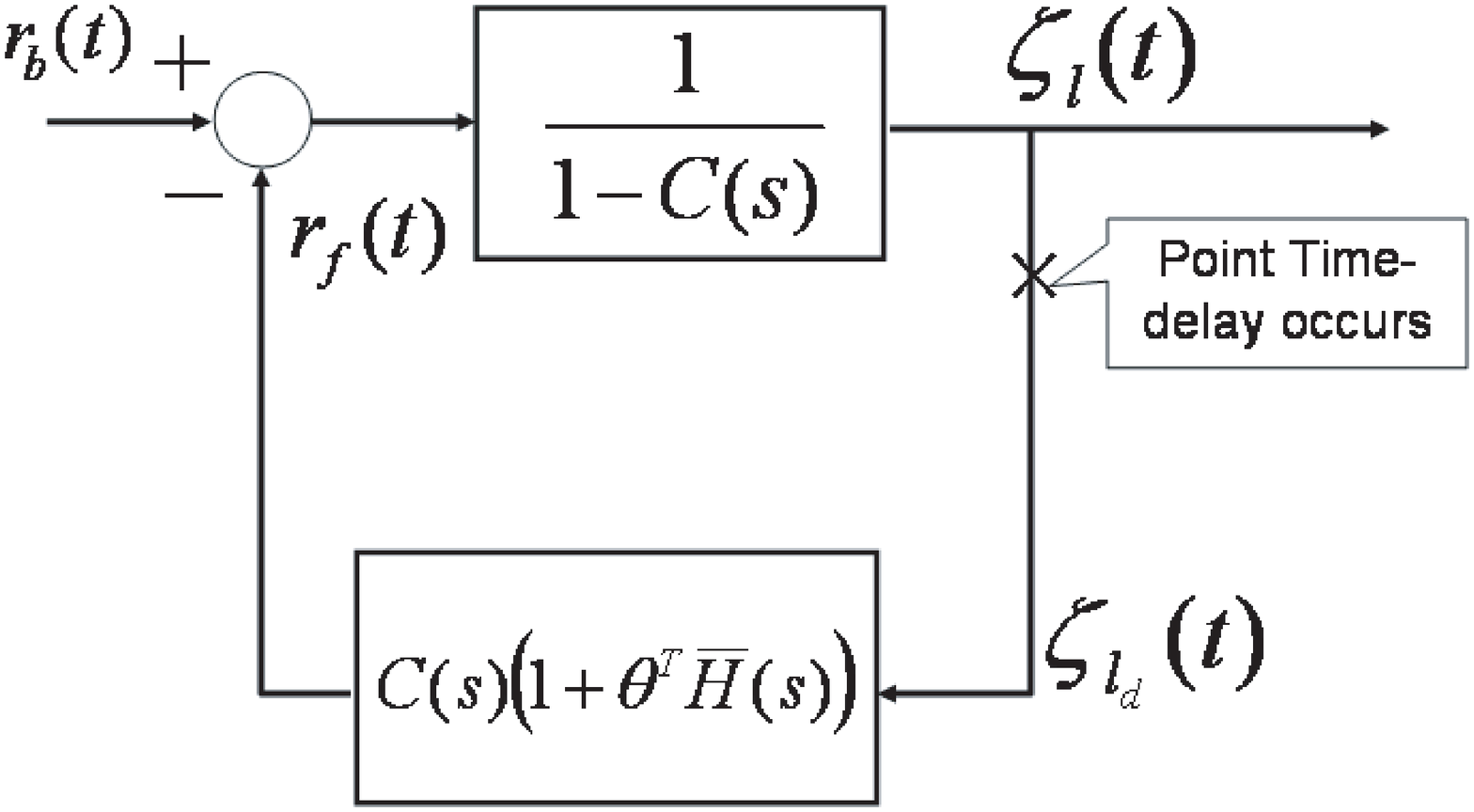}
\caption{LTI system}\label{fig:newLTIsys}
\end{center}
\end{figure}

We notice that
 the phase margin of this LTI system
can be determined by its open-loop transfer function from
$\zeta_{l_d}(t)$ to $\zeta_l(t)$. It follows from
(\ref{LTImarginori3}), (\ref{LTImarginori2}), and
(\ref{LTImarginori1}) that $ \omega u_l(s) =\frac{C(s)\left( k_g
r(s)-\zeta_{l_d}(s) -\theta^{\top} \bar{H}(s)
\zeta_{l_d}(s)\right)-\omega \epsilon_l(s)}{1-C(s)}\,, \nonumber $
and hence (\ref{LTImarginori4}) implies that $ \zeta_l(s)   =
\frac{C(s)\left( k_g r(s)-\zeta_{l_d}(s) -\theta^{\top} \bar{H}(s)
\zeta_{l_d}(s)\right)-\omega \epsilon_l(s)}{1-C(s)} +\sigma(s)$.
Therefore, it can be equivalently written as:
\begin{eqnarray}
\zeta_l(s) & = & \frac{1}{1-C(s)} \left(r_b(s)-
r_f(s)\right)\,,\nonumber\\
r_f(s) & = &  C(s)(1+\theta^{\top}\bar{H}(s)) \zeta_{l_d}(s)\,,\label{newLTIsys}\\
r_b(s) & = & C(s) k_g r(s)+(1-C(s))\sigma(s)-\omega
\epsilon_l(s)\,.\nonumber
\end{eqnarray}
Assume that $\tilde{r}_l(t)$ is such that $\epsilon_l(t)$ is
bounded. Since $\sigma(t)$ and $r(t)$ are  bounded, $C(s)$ is
strictly proper and stable,  then $r_b(t)$ is also bounded. The
block-diagram of the closed-loop system in (\ref{newLTIsys}) is
shown in Figure \ref{fig:newLTIsys}.

The open-loop transfer function of the system in (\ref{newLTIsys})
is:
\begin{equation}\label{Hosdefinion}
H_o(s) =  C(s)(1+\theta^{\top}\bar{H}(s)) /(1-C(s))\,,
\end{equation}
the phase margin $\mathcal{P}(H_o(s))$ of which can be derived
from its Bode plot easily.  Its time-delay margin  is given by:
\begin{equation}\label{TmarginDef}
\mathcal{T}(H_o(s))=\mathcal{P}(H_o(s))/\omega_c\,,
\end{equation}
where $\mathcal{P}(H_o(s))$ is the phase margin of the open-loop
system $H_o(s)$, and $\omega_c$ is the cross-over frequency of
$H_o(s)$. The next lemma states sufficient condition for
boundedness of all the states in the system
(\ref{LTImarginori3})-(\ref{LTImarginori4}), including the
internal states.
\begin{lem}\label{lem:LTItmori}
Let
\begin{equation}\label{thm_new000}
\tau < \mathcal{T}(H_o(s))\,
\end{equation}
and $\epsilon_b$ be any
positive number such that $\Vert
\epsilon_l\Vert_{\mathcal{L}_{\infty}}\leq \epsilon_b$. Then the
signals $ \zeta_l(t)$, $ x_l(t)$, $ u_l(t)$,  $\eta_l(t)$  are
bounded.
\end{lem}
{\bf Proof:} Since $\epsilon_l(t)$ is bounded and $\tau  <
\mathcal{T}(H_o(s))$, then boundedness of $\zeta_l(t)$ follows
from definition of $\mathcal{T}(H_o(s))$. Boundedness of
$\zeta_{l_d}(t)$ follows from its definition in (\ref{zeroinit4}).
Since $\zeta_l(t)$ and $\sigma(t)$  are bounded, it follows from
(\ref{LTImarginori4}) that $ u_l(t)$ is bounded, and
(\ref{LTImarginori1}) implies boundedness of
 $\eta_l(t)$. Notice that since $u_l(t)$ and $\epsilon_l(t)$ are
bounded, it follows from (\ref{LTImarginori2}) that $\theta^{\top}
x_l(t)$ is bounded. We notice that $x_l(s)$ in
(\ref{LTImarginori3}) can be written as $ x_l(s)=H(s)
(\theta^{\top} x_l(s)+\zeta_{l_d}(s)) $, which leads to
boundedness of $x_l(t)$.  $\hfill{\square}$

For any $\tau< \mathcal{T}(H_o(s))$ and any $\epsilon_b>0$, Lemma
\ref{lem:LTItmori} guarantees that the  map $\Delta_{n}:
{\rr}^{+}\times [0,{\mathcal{T}}(H_o(s))) \rightarrow {\rr}^{+}$
\begin{eqnarray}
\Delta_{n}(\epsilon_b,\tau) & = & \max_{\Vert
\epsilon_l\Vert_{\mathcal{L}_{\infty}} \leq \epsilon_b} \Vert
\sigma+ \eta_l\Vert_{\mathcal{L}_{\infty}}\label{Deltandef}
\end{eqnarray}
is well defined. We note that strictly speaking $\eta_l(t)$
depends not only on $\epsilon_l(t)$ and $\tau$, but also upon
other arguments, like $\sigma(t)$ and other variables of the
system that are used for definition of $\eta_l(t)$. These are
dropped due to their non-crucial role in the subsequent analysis.

\begin{lem}\label{lem:LTItm}
Let $\tau$ comply with (\ref{thm_new000}), and $\epsilon_b$ be any
positive number. If $\tilde{r}_l(t)$ is such that the resulting $
\epsilon_l (t)$ is bounded
\begin{equation}
 \Vert \epsilon_l \Vert_{\mathcal{L}_{\infty}}  \leq
\epsilon_b\,,\label{lem_LTItm000}
\end{equation}
and
\begin{equation}
 2\omega \Vert u_l
\Vert_{\mathcal{L}_{\infty}}+2 L \Vert x_l
\Vert_{\mathcal{L}_{\infty}}+ 2 \Delta \geq
 \Vert \tilde{r}_l
\Vert_{\mathcal{L}_{\infty}}\,,\label{lem_LTItm0}
\end{equation}
where
\begin{equation}
\Delta  =  \Delta_n(\epsilon_b,\tau)+\delta_1 \,,\label{thm_222}
\end{equation}
$\delta_1$ is arbitrary positive constant,
 then $\eta_l(t)$ is differentiable and the
$\mathcal{L}_{\infty}$ norm of $\dot{\eta}_l(t)$ is finite.
\end{lem}
{\bf Proof:} It follows from (\ref{lem_LTItm000}) and Lemma
\ref{lem:LTItmori} that  $x_l(t)$, $u_l(t)$,
$\Delta_n(\epsilon_b,\tau)$ are bounded. Hence, it follows from
(\ref{lem_LTItm0}) that $\tilde{r}_l(t)$ is also bounded. Since
$C(s)$ is strictly proper and stable, bounded $\tilde{r}_l(t)$
ensures that $\epsilon_l(t)$ is differentiable with bounded
derivative. Using similar methods, we prove that both $u_l(t)$ and
$\zeta_{l_d}(t)$ have bounded derivatives. Since $\dot{\sigma}(t)$
is bounded, it follows from (\ref{LTImarginori1}) that $\dot
\eta_l(t)$ is bounded.
 $\hfill{\square}$

For any $\tau< \mathcal{T}(H_o(s))$ and any $\epsilon_b>0$, Lemma
\ref{lem:LTItm} guarantees that the following map $\Delta_{d}:
{\rr}^{+}\times [0,{\mathcal{T}}(H_o(s))) \rightarrow {\rr}^{+}$
\begin{eqnarray}
\Delta_{d}(\epsilon_b,\tau)= \max_{\tilde r_l(t)}\Vert
\dot{\sigma}+\dot{\eta}_l\Vert_{\mathcal{L}_{\infty}}\label{Deltaddef}
\end{eqnarray}
is well defined, where $\tilde{r}_l(t)$  complies with
(\ref{lem_LTItm000}) and (\ref{lem_LTItm0}). Further, let
\begin{eqnarray}
\theta_{m}(\epsilon_b,\tau)  & \triangleq & \max_{\theta\in
\Theta} \sum_{i=1}^{n} 4 \theta_i^2+ 4\Delta^2+4
\left(\omega_{u}-\omega_{l} \right)^2
\nonumber\\
& & +2\lambda_{\max}(P)
\Delta_{d}(\epsilon_b,\tau)\Delta/\lambda_{\min}(Q)
\,,\label{thetamaxdef2} \\
\epsilon_c(\epsilon_b,\tau) & = &\Big\| C(s) \frac{1}{c_{o}^{\top}
H(s)} c_{o}^{\top} \Big\|_{\mathcal{L}_1}
\sqrt{\frac{\theta_{m}(\epsilon_b,\tau)}{\lambda_{\max}(P)
\epsilon_b^2}}\,.\label{epsilon_cdef}
\end{eqnarray}
 We notice that for any finite $\epsilon_b\in \rr^{+}$
and any $\tau$ verifying (\ref{thm_new000}), we have finite
$\Delta_{n}(\epsilon_b,\tau)$ and $\Delta_{d}(\epsilon_b,\tau)$,
and hence finite $\epsilon_c(\epsilon_b,\tau)$, if
$\tilde{r}_l(t)$ complies with (\ref{lem_LTItm000}) and
(\ref{lem_LTItm0}).

\subsection{Time-delay Margin of the Closed-loop Adaptive System}

In this section we analyze the time-delay margin  for the
closed-loop adaptive system with the $\mathcal{L}_1$ adaptive
controller.
The main result is given by the following theorem.
\begin{thm}\label{thm:phaseeq}
Consider the closed-loop adaptive system, comprised of System 1 in
(\ref{L1_td1})-(\ref{zeroinit1})  with (\ref{L1_companionmodal}),
(\ref{adaptivelaw_L1})-(\ref{adaptivelaw_L3}), (\ref{controllaw})
and the LTI system in (\ref{LTImarginori3})-(\ref{LTImarginori4})
in the presence of the same time delay $\tau$. For any
$\epsilon_b\in \rr^{+}$ choose the set $\Delta$ as in
(\ref{thm_222}) and
\begin{equation}
\Gamma_c \geq
\sqrt{\epsilon_c(\epsilon_b,\tau)}+\delta_2\,,\label{thm_333}
\end{equation}
where $\delta_2$ is arbitrary positive constant. Then for every
$\tau$ satisfying $ \tau < \mathcal{T}(H_o(s)) $,
 there
exists  exogenous signal $\tilde{r}_l(t)$ ensuring $ \Vert
\epsilon_l\Vert_{\mathcal{L}_\infty} < \epsilon_b\,, $ and
\begin{equation}\label{thm_con}
x_l(t) =x_d(t)\,,\quad  u_l(t)= u(t)\,,\qquad \forall t\geq 0\,.
\end{equation}
\end{thm}
%
The proof  of this Theorem  is given in the Appendix.
Theorem \ref{thm:phaseeq} establishes the  equivalence of state
and control trajectories of the closed-loop adaptive system and
the
 LTI system in (\ref{LTImarginori3})-(\ref{LTImarginori4}) in
the presence of the same time-delay. Therefore the time-delay
margin of the system in
(\ref{LTImarginori3})-(\ref{LTImarginori4}) can be used as a
conservative lower bound for the time-delay margin of the
closed-loop adaptive system.

\begin{cor}\label{cor:td}
Given the system in (\ref{problemnew}) and the $\mathcal{L}_1$
adaptive controller defined via (\ref{L1_companionmodal}),
(\ref{adaptivelaw_L1})-(\ref{adaptivelaw_L3}) and
(\ref{controllaw}) subject to (\ref{condition3}), where $\Gamma_c$
and $\Delta$ are large enough, the closed-loop adaptive system is
stable in the presence of time delay $\tau$ in  its output if
$
\tau  < \mathcal{T}(H_o(s))\,,
$
where $\mathcal{T}(H_o(s))$ is defined in (\ref{TmarginDef}).
\end{cor}
\indent{} The proof of Corollary \ref{cor:td} follows from Lemma
\ref{lem:LTItmori} and Theorem \ref{thm:phaseeq} directly.


\section{Gain Margin Analysis}\label{sec:gainmargin}
We now analyze the gain margin of the system in (\ref{problemnew})
with $\mathcal{L}_1$ adaptive controller. By inserting a gain
module $g$ into the control loop, the system in (\ref{problemnew})
can be formulated as:
\begin{equation}
\dot{x}(t)  =   A_m x(t)+b \left(  \omega_g u(t)+\theta^{\top}(t)
x(t)+ \sigma(t) \right)\,, \label{gm_2}
\end{equation}
where  $\omega_g=g \omega$. We note that this transformation
implies that the set $\Omega$ in the application of the Projection
operator for adaptive laws needs to increase accordingly. However,
increased $\Omega$ will not violate the stability criterion. Thus,
it follows from (\ref{eqn:Delta}) that the gain margin of the
$\mathcal{L}_1$ adaptive controller is determined by:
\begin{equation}\label{Gmdef}
\mathcal{G}_m =[\omega_l/\omega_{l_0}, \;\;\omega_u/\omega_{u_0}].
\end{equation}
If $ g \in \mathcal{G}_m\,, $ then the closed-loop system in
(\ref{gm_2}) satisfies the $\mathcal{L}_1$ stability criterion,
implying that the entire closed-loop system is stable. We note
that the lower-bound of $\mathcal{G}_m$ is greater than zero. Eq.
(\ref{Gmdef}) implies that arbitrary gain margin can be obtained
through appropriate choice of $\Omega$.
%

\section{Main Results}\label{sec:generalization}

Combining the results of Theorem \ref{thm:5}, and Theorems of
Sections \ref{sec:timedelaymargin} and \ref{sec:gainmargin}, we
have the following results:

\begin{thm}\label{thm:main}
Given the system in (\ref{problemnew}) and the $\mathcal{L}_1$
adaptive controller defined via (\ref{L1_companionmodal}),
(\ref{adaptivelaw_L1})-(\ref{adaptivelaw_L3}) and
(\ref{controllaw}) subject to (\ref{condition3}), we have:
\begin{eqnarray}
&&\lim_{\Gamma_c\rightarrow \infty} \left(x(t)-x_{ref}(t) \right)
=
 0\, ,\qquad \forall t\geq 0,\label{cor1_00}
\\
&&\lim_{\Gamma_c\rightarrow \infty} \left(u(t)-u_{ref}(t)\right)
=
 0\,,\qquad \forall t\geq 0\,,\label{cor1_01}\\
&&\lim_{\Gamma_c\rightarrow \infty} \mathcal{T}  \geq
\mathcal{T}(H_o(s))  \,,\quad
 \mathcal{G}  \supseteq  \mathcal{G}_m
\,,\label{cor1_04}
\end{eqnarray}
where $\mathcal{T}$ and $\mathcal{G}$ are the time-delay and gain
margins of the $\mathcal{L}_1$ adaptive controller,  while
$\mathcal{T}(H_o(s))$, $\mathcal{G}_m$ are defined in
(\ref{TmarginDef}) and (\ref{Gmdef}).
\end{thm}
\indent{  } The inequalities in  (\ref{cor1_04}) imply that
$\mathcal{T}(H_o(s))$ and $\mathcal{G}_m$ are just conservative
bounds of the stability margins.

\section{Simulations}\label{sec:simu}
We consider the same system from \cite{CDC06_chengyu_L1}, in which
a single-link robot arm is rotating on a vertical plane. Assuming
constant $\theta(t)$,  it can be cast into the form in
(\ref{problemnew}) with $ A_m= \left[
\begin{array}{cc}
0 &  1\\
-1 & -1.4
\end{array}
\right]\,,\quad b=\left[
\begin{array}{c}
0 \\
1
\end{array}
\right]\,, \quad c=\left[
\begin{array}{c}
1 \\
0
\end{array}
\right]$. Let $ \theta=[2\;\;2]^{\top},\,\,
\omega=1,\,\,\sigma(t)=\sin(\pi t)\,, $ so that the compact sets
can be conservatively chosen as $ \Omega_0=[0.2,~\,5],\,
\Theta=[-10,\,10],\, \Delta_0=[-10,\,10]\,, $ respectively. Next,
we analyze the stability margins of the $\mathcal{L}_1$ adaptive
controller for this system numerically.
\begin{figure}[!h]
\begin{center}
\includegraphics[width=2.5in,height=1.8in]{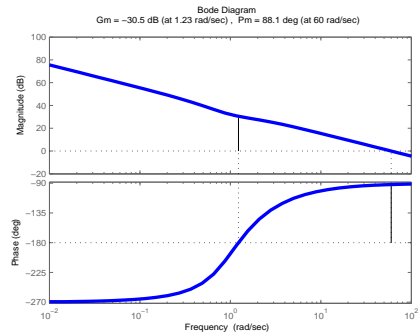}
\caption{Bode plot of $H_o(s)$ for $\theta=[2\;\;2]^{\top}$,
$\omega=1$}\label{fig:phasebode}
\end{center}
\end{figure}

\begin{figure}[h!]
\centering \mbox{\subfigure[ $x_1(t)$ (solid),  $\hat{x}_1(t)$
(dashed), and $r(t)$(dotted)
 ]{\epsfig{file=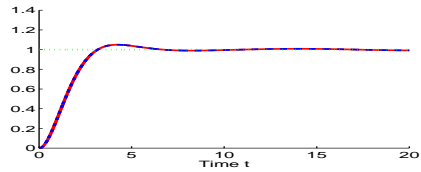,width=2.5in,height=0.9in
}\label{fig:L1td_y}}} \mbox{ \subfigure[ Time-history of $u(t)$
]{\epsfig{file=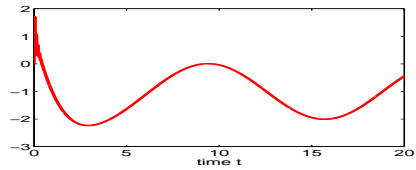,width=2.5in,
height=0.9in}\label{fig:L1td_u} }} \caption{Performance of
$\mathcal{L}_1$ adaptive controller with time-delay $0.02s$}
\end{figure}

For $ \theta=[2\,\,\,\,2]^{\top},~ \omega=1 $ we can derive
$H_o(s)$ in (\ref{Hosdefinion})  and look at its Bode plot in Fig.
\ref{fig:phasebode}. It has phase margin $88.1^{\circ}(1.54
\rm{rad})$ at cross frequency $9.55 {\rm Hz}(60\rm{rad/s})$.
Hence, the time-delay margin can be derived from
(\ref{TmarginDef}) as: $ \mathcal{T}(H_o(s))
=\frac{1.54\rm{rad}}{60\rm{rad/s}}=0.0256$. We set $
\Delta=[-1000\;\;\; 1000]^{\top}, \quad \Gamma_c=500000\,, $ and
run the $\mathcal{L}_1$ adaptive controller with time-delay
$\tau=0.02$. The simulations  in Figs.
\ref{fig:L1td_y}-\ref{fig:L1td_u} verify  Corollary \ref{cor:td}.
As stated in Theorem \ref{thm:main}, the time-delay margin of the
LTI system  in (\ref{Hosdefinion}) provides only a conservative
lower bound for the time-delay margin of the closed-loop adaptive
system. So, we simulate the $\mathcal L_1$ adaptive controller in
the presence of larger time-delay, like $\tau=0.1$ sec., and
observe that the system is not losing its stability. Since
$\theta$ and $\omega$ are unknown to the controller, we derive the
$\mathcal{T}(H_o(s))$ for all possible $\theta\in\Theta$ and
$\omega\in \Omega$ and use the most conservative value. It gives $
\mathcal{T}(H_o(s)) = 0.005 s$. The gain margin can be arbitrarily
large as stated in (\ref{cor1_04}).

\section{Conclusion}\label{sec:conclusion}
In this paper, we derive the stability margins of $\mathcal L_1$
adaptive controller presented in \cite{CDC06_chengyu_L1}. To the
best of our knowledge, this is the first attempt to quantify the
time-delay margin for general closed-loop adaptive systems. With
the particular architecture presented in this paper, we prove that
increasing the adaptive gain leads to improved transient tracking
with improved stability margins. This presents a significant
improvement over conventional adaptive control schemes, in which
increasing the adaptive gain leads to reduced tolerance to
time-delay in input/output channels.

\section*{Appendix}
{\bf Proof of Theorem \ref{thm:phaseeq}:} In the closed-loop
adaptive system in (\ref{L1_tv})  for any $t\geq 0$, we notice
that if $ \Vert (\sigma+\eta)_t\Vert_{\mathcal{L}_{\infty}} \leq
\Delta\,, $ and $\dot{\sigma}(t),~\dot{\eta}(t)$  have finite
derivatives over $[0, ~t]$, then application of $\mathcal{L}_1$
adaptive controller from \cite{CDC06_chengyu_L1} is well-defined.
Let $d_{t}$ denote the truncated $\mathcal{L}_{\infty}$ norm
\begin{eqnarray}\label{dtdef}
d_{t} = \Vert
(\dot{\sigma}+\dot{\eta})_{t}\Vert_{\mathcal{L}_{\infty}}\,.
\end{eqnarray}
It follows from (\ref{L1_companionmodal}) and (\ref{L1_tv}) that
$\tilde x_q(s)=H(s) \tilde r(s)$, where $\tilde x_q(s)$ and
$\tilde r(s)$ are the Laplace transformations of $\tilde x_q(t)
=\hat{x}(t)-x_q(t)$ and
\begin{equation}\label{mainrtdef}
\tilde{r}(t)=\tilde{\omega}(t) u_q(t) +\tilde{\theta}^{\top}(t)
x_q(t)+ \tilde{\sigma}(t)\,.
\end{equation}
This along with Eq.  (50)
 in \cite{CDC06_chengyu_L1} implies that
\begin{eqnarray}\label{L1mu}
u_q(t)_{t\in[0,t]} & = & \mathcal{S}(\frac{C(s)}{\omega}, u_q(0),
(k_g r(t)-\theta^{\top}
x_q(t)-\nonumber\\
& &
\sigma(t)-\eta(t)-\tilde{r}(t))_{t\in[0,t)})\,,\nonumber\\
\tilde{x}_q(t)_{t\in[0,t]} & = & \mathcal{S}(H(s), \tilde{x}_q(0),
\tilde{r}(t)_{t\in[0,t)})\,,
\end{eqnarray}
where
$\tilde{\sigma}(t)=\hat{\sigma}(t)-(\sigma(t)+\eta(t))$. Equation
(\ref{L1mu}) implies that
\begin{eqnarray}\label{nn_43}
 u_q(t)_{t\in[0,t]}  & = &  \mathcal{S}(\frac{C(s)}{\omega},  u_q(0),
(k_g r(t)-\theta^{\top}
x_q(t)-\nonumber\\
& & \sigma(t)-\eta(t))_{t\in[0,t)})-\epsilon(t)_{t\in[0,t]}\,,
\end{eqnarray}
where
\begin{equation}\label{nn_44}
\epsilon(t)_{t\in[0,t]}= \mathcal{S}(C(s)/\omega, 0,
\tilde{r}(t)_{t\in[0,t)})\,.
\end{equation}
 We
further define
\begin{equation}
 \theta_{t}  \triangleq \max_{\theta\in \Theta} \sum_{i=1}^{n} 4
\theta_i^2+\max_{\sigma\in \Delta} 4 \sigma^2 +4
\left(\omega_u-\omega_l
\right)^2+2\frac{\lambda_{\max}(P)}{\lambda_{\min}(Q)}
 d_{t}  \Delta
 \,,\label{thetamaxdeftv}
\end{equation}
where $d_{t}$ is defined in (\ref{dtdef}). It can be verified
easily that Lemma \ref{lem:1}  holds for truncated norms as well
so that $ \Vert \tilde{x}_{q_t} \Vert_{\mathcal{L}_{\infty}} \leq
\sqrt{\frac{\theta_{t}}{\lambda_{\min}(P) \Gamma_c}}$. Since $
\epsilon(s) =  \frac{C(s)}{\omega c_{o}^{\top} H(s)} c_{o}^{\top}
H(s) \tilde{r}(s) =  \frac{C(s)}{\omega c_{o}^{\top} H(s)}
c_{o}^{\top} \tilde{x}_q(s)\,, $ then $\epsilon(t)$ can be upper
bounded as
\begin{equation}\label{L1tilder}
\Vert \epsilon_{t} \Vert_{\mathcal{L}_{\infty}} \leq \Big\| C(s)
\frac{1}{\omega c_{o}^{\top} H(s)} c_{o}^{\top}
\Big\|_{\mathcal{L}_1} \sqrt{\frac{\theta_{t}}{\lambda_{\min}(P)
\Gamma_c}}\,.
\end{equation}

In the three steps below, we prove the existence of a continuously
differentiable $\eta(t)$ with uniformly bounded derivative in the
closed-loop adaptive system (\ref{L1_tv}),
(\ref{L1_companionmodal}),
(\ref{adaptivelaw_L1})-(\ref{adaptivelaw_L3}), (\ref{controllaw})
and the existence of $r_l(t)$ in time-delayed LTI system  such
that for any $t\geq 0$,
\begin{eqnarray}
|\sigma(t)+\eta(t)|  <  \Delta\,, x_{o}(t)  =  x_q(t) \,, \qquad \qquad \qquad &&\label{disinseti}\\
 \Vert \epsilon_{l_t}\Vert_{\mathcal{L}_{\infty}} <
\epsilon_b\,,x_l(t) =x_q(t),u_l(t)=u_q(t),
\epsilon_l(t)=\epsilon(t)\,. &&\label{step0_0}
\end{eqnarray}
With (\ref{disinseti}), Lemma \ref{lem:equivalence} implies that $
x_{d}(t)  =  x_q(t),u(t) =u_q(t)$ for any $t\geq 0$, which
combining (\ref{step0_0}) proves Theorem \ref{thm:phaseeq}.

\underline{\it Step 1:}  Let
\begin{equation}\label{step1_88}
\zeta(t)= \omega u_q(t)+\sigma(t)\,.
\end{equation}
We further define
\begin{eqnarray}
& & \zeta_d =\left\{
\begin{array}{cc}
0 \,,& t\in[0, \tau)\\
\zeta(t-\tau) \,,& t\geq \tau
\end{array}
\right.\,.\label{thm_step1_00}
\end{eqnarray}
Since (\ref{zeroinit0}) and (\ref{zeroinit1}) imply that $
x_{o}(t)=0$ for any $t\in [0,\tau]$, it follows from
(\ref{thm_step1_00}) and the definition of the map $\mathcal{S}$
that $ x_{o}(t)|_{t\in[0 ,\tau]}  =  \mathcal{S}\left( \bar{H}(s),
x_{o}(0),
 \zeta_{d}(t)_{t\in[0,\tau)}\right).
$ For $i\geq 1$, it follows from the definition of the
time-delayed open-loop system that
\begin{equation}
  x_{o}(t)|_{t\in[i \tau,(i+1)\tau]}  =  \mathcal{S}\left(
\bar{H}(s), x_{o}(i \tau),
 \zeta_{d}(t)_{t\in[i \tau,(i+1)\tau)}\right)\,.\label{step2_14}
\end{equation}
Hence, (\ref{step2_14}) holds for any $i$. We note that
(\ref{nn_44}) implies that $ \epsilon(0)=0\,. $
These along with (\ref{zeroinit0}), (\ref{zeroinit1}),
(\ref{zeroinit4}), (\ref{zeroinit3new}), imply that for $i=0$
\begin{eqnarray}
&  & u_q(i \tau)  =   u_l(i \tau),
 \epsilon(i \tau)   =  \epsilon_l(i \tau),
 x_{o}(i \tau)  =  x_q(i \tau) = x_l(i \tau)\,,\nonumber\\
& &   \zeta_d(t)  =  \zeta_{l_d}(t)\,,t< (i+1)
 \tau\,,\;\;
 |\epsilon(t) |  <  \epsilon_b\,,\quad t\leq i
 \tau\,.\nonumber
\end{eqnarray}

\underline{\it Step 2:} Assume that for any $i$ the following
conditions hold:
\begin{eqnarray}
&&  u_q(t)  =   u_l(t)\,,\quad t\leq i \tau\,,\label{step2_11} \\
&& \epsilon(t)  =  \epsilon_l(t)\,,\quad t= i \tau\,,\label{step2_18} \\
&& x_{o}(t)  =  x_q(t) = x_l(t)\,,\quad t \leq i \tau\,,\label{step2_12}\\
&& \zeta_d(t)  =  \zeta_{ld}(t)\,,~ \forall ~t\in [i\tau, (i+1)
\tau)\,,\label{step2_13}\\
 & & |\epsilon(t)|  <  \epsilon_b\,,~ \forall ~ t\leq i \tau\,.\label{step2_16}
\end{eqnarray}
For $i\ge 1$, further  assume that there
 exist
bounded $\tilde{r}_l(t)$ and continuously differentiable $\eta(t)$
with bounded derivative over $t\in[0, i\tau)$ such that $\forall
~t< i \tau$
\begin{eqnarray}
 \eta(t) = \eta_l(t)\,, \quad |\sigma(t)+\eta(t)| <   \Delta\,. \label{step2_19}
\end{eqnarray}
We prove below that there exist bounded $\tilde{r}_l(t)$ and
continuously differentiable $\eta(t)$ with bounded derivative over
$t\in[0, (i+1)\tau)$ such that (\ref{step2_11})-(\ref{step2_19})
hold for $i+1$, too.

We note that (\ref{LTImarginori3}) implies that
\begin{equation}\label{naira1}
x_{l}(t)|_{t\in[i\tau,(i+1)\tau]}  =  \mathcal{S}\left(
\bar{H}(s), x_l(i \tau),
 \zeta_{l_d}(t)_{t\in[i \tau,(i+1)\tau)}\right).
\end{equation}
Using (\ref{step2_12})-(\ref{step2_13}), it follows from
(\ref{step2_14}) and (\ref{naira1})  that
\begin{equation}\label{step2_33}
x_{o}(t)=x_l(t),\quad \forall~ t\in[i \tau,\,(i+1)\tau]\,.
\end{equation}
We assumed in (\ref{step2_19}) that if $i\geq 1$, then there
exists $\eta(t)$ over $[0,i\tau)$. We now define $\eta(t)$ over $
[i \tau, (i+1) \tau)$ as:
\begin{equation}\label{step2_6}
\eta(t) = \zeta_d(t)-\omega u_q(t)-\sigma(t)\,,~~t\in[i \tau,
(i+1) \tau)\,.
\end{equation}
Since (\ref{L1_tv}) implies that $
x_{q}(t)|_{t\in[i\tau,(i+1)\tau]}   =   \mathcal{S}( \bar{H}(s),
x_q(i \tau) (\omega u_q(t)+\sigma(t)+\eta(t))_{t\in[i
\tau,(i+1)\tau)})$, it follows from (\ref{step2_6}) that $
x_{q}(t)|_{t\in[i\tau,(i+1)\tau]}   =   \mathcal{S}( \bar{H}(s),
x_q(i \tau),\zeta_d(t)_{t\in[i \tau,(i+1)\tau)})$. Along with
(\ref{step2_14}) and (\ref{step2_12})
 this ensures that
\begin{equation}\label{xqdeq}
x_q(t) = x_{o}(t),~~ \forall ~t\in[i\tau, (i+1) \tau]\,.
\end{equation}
However, the definition in (\ref{step2_6}) does not guarantee
\begin{equation}\label{step2_666}
|\sigma(t)+\eta(t)| < \Delta\,,~~t\in[i \tau, (i+1) \tau)\,,
\end{equation}
which is required for application of $\mathcal L_1$ adaptive
controller.

We prove (\ref{step2_666}) by contradiction. Since $\eta(t)$ is
continuous over $[i \tau, (i+1) \tau)$, if (\ref{step2_666}) is
not true, there must exist $t'\in [i \tau, (i+1) \tau)$ such that
\begin{eqnarray}
|\sigma(t)+\eta(t)| & < & \Delta \,,\quad \forall t< t'\,,\label{step2_n11}\\
|\sigma(t')+\eta(t')| & = & \Delta \,.\label{step2_n12}
\end{eqnarray}
It follows from (\ref{step2_14}) and (\ref{step2_6}) that $
x_{o}(t)_{t\in[i\tau, t']} = \mathcal{S}\Big( \bar{H}(s), x_{o}(i
\tau), (\omega u_q(t) +\sigma(t)+\eta(t))_{t\in[i \tau,
t')}\Big)$. It follows from (\ref{nn_43}) and (\ref{nn_44}) that
\begin{eqnarray}
& &  u_q(t)|_{t\in[i\tau,\,t']} =
\mathcal{S}\Big(\frac{C(s)}{\omega},u_q(i \tau)+\epsilon(i
\tau),(k_g r(t)-\theta^{\top}
x_q(t) \nonumber\\
& & -\sigma(t)-\eta(t))_{t\in[i\tau,\,t')}\Big)
-\epsilon(t)|_{t\in[i\tau,\,t']} \,,\label{step2_mu}
\end{eqnarray}
where
\begin{equation}\label{step2_new88}
\epsilon(t)|_{t\in[i\tau,\,t']}=\mathcal{S}\Big(C(s)/\omega,\epsilon(i\tau),
\tilde{r}(t)_{t\in[i\tau,\,t')}\Big)\,.
\end{equation}
We notice that if $i\geq 1$, then on $[0,i\tau)$ we have
$\tilde{r}_l(t)$ well defined. Let
\begin{equation}\label{step2_epsdef}
\tilde r_l(t) = \tilde{r}(t)\,, \quad t\in[i \tau,\,t')\,.
\end{equation}
We have $ \epsilon_l|_{t\in[i\tau,\,t']} =
\mathcal{S}\Big(C(s)/\omega,\epsilon_l(i\tau),
\tilde{r}(t)_{t\in[i\tau,\,t')}\Big)$, which along with
(\ref{step2_18}) and (\ref{step2_new88}) imply that
\begin{equation}\label{step2_nn22}
\epsilon_l(t)= \epsilon(t), ~~~\forall ~t\in[i \tau, t']\,.
\end{equation}
Hence,  (\ref{step2_11}), (\ref{step2_33}),
(\ref{xqdeq}), 
(\ref{step2_mu}) yield
\begin{eqnarray}
& & u_q(t)|_{t\in[i\tau,\,t']} = \mathcal{S}\Big(C(s)/\omega,
u_l(i \tau)+\epsilon(i \tau),  (k_g r(t) \nonumber\\
&& -\theta^{\top}
x_{l}(t)-\sigma(t)-\eta(t))_{t\in[i\tau,\,t')}\Big)-\epsilon(t)|_{t\in[i\tau,\,t']}\,.\label{step2_57}
\end{eqnarray}
It follows from (\ref{step2_nn22}) and (\ref{step2_57}) that
\begin{eqnarray}
& & u_q(t)|_{t\in[i\tau,\,t']} = \mathcal{S}\Big(C(s)/\omega,
u_l(i \tau)+\epsilon_l(i \tau),  (k_g r(t) \nonumber\\
&& -\theta^{\top}
x_{l}(t)-\sigma(t)-\eta(t))_{t\in[i\tau,\,t')}\Big)-\epsilon_l(t)|_{t\in[i\tau,\,t']}\,.\label{step2_56}
\end{eqnarray}
It follows from (\ref{LTImarginori1}) and (\ref{step2_13}) that
\begin{equation}\label{step2_55}
\eta_l(t) = \zeta_d(t)-\omega u_l(t)-\sigma(t)\,,~~t\in[ i
\tau,t']\,,
\end{equation}
which along with (\ref{LTImarginori2}) yields
\begin{eqnarray}
& &  u_l(t)|_{t\in[i\tau,\,t']} =
\mathcal{S}\Big(\frac{C(s)}{\omega},u_l(i \tau)+\epsilon_l(i
\tau),(k_g r(t) \nonumber\\
& & -\theta^{\top}
x_l(t)-\sigma(t)-\eta_l(t))_{t\in[i\tau,\,t')}\Big)
-\epsilon_l(t)|_{t\in[i\tau,\,t']} \,.\label{step2_mul}
\end{eqnarray}
From (\ref{step2_6}), (\ref{step2_56}), (\ref{step2_55}) and
(\ref{step2_mul}), we have
\begin{eqnarray}
 u_q(t) & = &  u_l(t),\quad \forall
t\in[i\tau,\,t']\,\label{step2_66}\\
\eta(t) & = & \eta_l(t),\quad \forall
t\in[i\tau,\,t')\,.\label{step2_67}
\end{eqnarray}
It follows from (\ref{step2_19}) and  (\ref{step2_67}) that
\begin{eqnarray}
\eta(t) & = & \eta_l(t),\quad \forall
t\in[0,\,t')\,.\label{step2_677}
\end{eqnarray}
We now prove  by contradiction that
\begin{equation}\label{step2_con1}
|\epsilon(t)| < \epsilon_b\,, \quad \forall~t~\in [i \tau, t']\,.
\end{equation}
 If (\ref{step2_con1}) is not true, then since $\epsilon(t)$ is continuous,
there exists some $\bar{t}\in[i\tau, t']$ such that
\begin{eqnarray}
|\epsilon(t)| & < & \epsilon_b\,,\qquad \forall ~t\in[i
\tau,\bar{t})\,,\label{step2_con21} \\
|\epsilon(\bar{t}) |& = & \epsilon_b\,. \label{step2_con22}
\end{eqnarray}
It follows from (\ref{step2_16}) that
\begin{equation}\label{main_222}
|\epsilon(t)| \leq \epsilon_b\,, \qquad \forall~ [0, \bar{t}]\,.
\end{equation}
It follows from (\ref{step2_11}), (\ref{step2_12}),
(\ref{step2_33}), (\ref{xqdeq}) and (\ref{step2_66}) that $
u_q(t)=u_l(t)\,,\quad x_q(t)=x_l(t)$ for any $ t\in[0,\bar{t}]$.
Therefore, (\ref{mainrtdef}) and (\ref{step2_epsdef}) imply
 that
$ \tilde{r}_l(t) = \tilde{\omega}(t)
u_l(t)+\tilde{\theta}^{\top}(t) x_l(t)+ \tilde{\sigma}(t) $, and
hence
\begin{equation}\label{main111}
\Vert \tilde{r}_{l_{\bar{t}}} \Vert_{\mathcal{L}_{\infty}}  \leq 2
\omega  \Vert u_{l_{\bar{t}}}\Vert_{\mathcal{L}_{\infty}}+ L \Vert
x_{l_{\bar{t}}}\Vert_{\mathcal{L}_{\infty}}+ 2 \Delta\,.
\end{equation}
From (\ref{main_222}) and (\ref{main111}), Lemmas
\ref{lem:LTItmori} and \ref{lem:LTItm} imply that $\eta_l(t)$ is
bounded and differentiable with bounded derivative. Further,
 it follows from (\ref{Deltandef}) and (\ref{Deltaddef})
that
\begin{eqnarray}
| \sigma(t)+\eta_l(t) | & \leq &
\Delta_n(\epsilon_b,\tau)\,,\qquad
\forall t\in [0,\bar{t}]\,, \nonumber\\
| \dot{\sigma}(t)+ \dot{\eta}_l(t) | & \leq &
\Delta_d(\epsilon_b,\tau)\,,\qquad \forall t\in [0,\bar{t}]\,.
\end{eqnarray}
Since (\ref{step2_677}) holds,  $\eta(t)$ is also bounded and
differentiable with bounded derivative over $[0, t')$ and further
\begin{eqnarray}
| \sigma(t)+\eta(t) | & \leq & \Delta_n(\epsilon_b,\tau)\,,\qquad
\forall t\in [0,\bar{t}]\,, \label{main_88}\\
| \dot{\sigma}(t) +\dot{\eta}(t) | & \leq &
\Delta_d(\epsilon_b,\tau)\,,\qquad \forall t\in [0,\bar{t}]\,.
\label{main_8}
\end{eqnarray}
It follows from (\ref{L1tilder}) that
\begin{equation}\label{L1tildera}
\Vert \epsilon_{\bar{t}} \Vert_{\mathcal{L}_{\infty}} \leq \Big\|
C(s) \frac{1}{\omega c_{o}^{\top} H(s)} c_{o}^{\top}
\Big\|_{\mathcal{L}_1}
\sqrt{\frac{\theta_{\bar{t}}}{\lambda_{\min}(P) \Gamma_c}}\,.
\end{equation}
It follows from (\ref{thetamaxdef2}), (\ref{thetamaxdeftv})  and
(\ref{main_8}) that
\begin{equation}\label{main_5}
\theta_{\bar{t}} \leq \theta_{m}(\epsilon_b,\tau)\,.
\end{equation}
Hence, (\ref{L1tilder}) and (\ref{main_5}) imply that $ \Vert
\epsilon_{\bar{t}} \Vert_{\mathcal{L}_{\infty}} \leq \Big\| C(s)
\frac{1}{\omega c_{o}^{\top} H(s)} c_{o}^{\top}
\Big\|_{\mathcal{L}_1}
\sqrt{\frac{\theta_{m}(\epsilon_b,\tau)}{\lambda_{\min}(P)
\Gamma_c}}$. From (\ref{epsilon_cdef}) and (\ref{thm_333}) we have
$ \Vert \epsilon_{\bar{t}} \Vert_{\mathcal{L}_{\infty}} <
\epsilon_b$, which contradicts (\ref{step2_con22}). Therefore,
(\ref{step2_con1}) holds.

If (\ref{step2_con1}) is true, it follows from (\ref{step2_16})
that
$$
|\epsilon(t)| < \epsilon_b\,,~~ \forall ~t\in[0, t']\,.
$$
Hence, it follows from (\ref{Deltandef}) and (\ref{step2_677})
that
\begin{equation}
| \sigma(t)+\eta(t) | \leq \Delta_n < \Delta\,,
\end{equation}
which contradicts (\ref{step2_n12}). Hence, we have
\begin{equation}\label{main_yy}
| \sigma(t)+\eta(t) | < \Delta, ~~\forall~ t\in[i \tau, (i+1)
\tau].
\end{equation}

Therefore, combining (\ref{step2_33}), (\ref{xqdeq}),
(\ref{step2_nn22}), (\ref{step2_66}), (\ref{step2_67}),
(\ref{step2_con1}), (\ref{main_yy}), we proved that there exist
$\tilde{r}_l(t)$ and continuously differentiable $\eta(t)$ in $[0,
(i+1) \tau)$, which ensures
\begin{eqnarray}
x_o(t) & = & x_q(t) = x_l(t), \forall~ t\in[i\tau, (i+1) \tau], \label{res1}\\
\epsilon(t) & = & \epsilon_l(t), ~~\forall~ t\in[i\tau, (i+1) \tau],\label{res2} \\
u_q(t) & = & u_l(t), ~~\forall~ t\in[i\tau, (i+1) \tau], \label{res3}\\
\eta(t) & = & \eta_l(t), ~~\forall~ t\in[i\tau, (i+1) \tau), \label{res4}\\
|\epsilon(t) | & < & \epsilon_b,~~\forall~t\in[0, (i+1) \tau], \label{res5}\\
| \sigma(t)+\eta(t) | & < & \Delta, ~~\forall~ t\in[0, (i+1)
\tau].\label{res6}
\end{eqnarray}

It follows from (\ref{LTImarginori4}), (\ref{step1_88})  and
(\ref{res3}) that
$$
\zeta(t) = \zeta_l(t) , \quad \forall ~t\in[i\tau, (i+1) \tau)\,.
$$
Therefore (\ref{zeroinit4}) and (\ref{thm_step1_00})   imply that
\begin{equation}\label{step2_kappa}
\zeta_d(t) = \zeta_{l_d}(t) , \quad \forall ~t\in[(i+1)\tau, (i+2)
\tau)\,.
\end{equation}
We note that Step 2 is proved in (\ref{res1})-(\ref{step2_kappa})
for $i+1$.

\underline{\it Step 3:} Step $1$ implies that the relationships
(\ref{step2_11})-(\ref{step2_16}) hold for $i=0$. By iterating the
results from Step $2$, we prove (\ref{disinseti})-(\ref{step0_0}),
which conclude proof of  the  Theorem. $\hfill{\square}$

\end{document}